\documentclass[reqno,A4paper]{amsart}

\usepackage{amsmath}
\usepackage{amssymb}
\usepackage{amsthm}
\usepackage{algorithm}
\usepackage{caption}
\usepackage{tikz} 
\usepackage{here}
\usepackage{hyperref}
\usepackage{algorithmic}
\usepackage{enumerate}
\usepackage{enumitem}
\usepackage{mathrsfs} 
\usepackage{eqlist}
\usepackage{array}

\theoremstyle{plain}
\newtheorem{thm}{Theorem}[section]
\newtheorem{pro}[thm]{Proposition}
\newtheorem{lem}{Lemma}[section]

\theoremstyle{definition}
\newtheorem{defn}{Definition}
\newtheorem{rem}{\textup{Remark}} 
\newtheorem{ex}{\textit{Example}} 

\numberwithin{equation}{section}

\newcommand{\CC}{\mathcal{C}}    
\newcommand{\A}{\mathcal{A}}    

\begin{document}

\title[The Comaximal filter graphs in residuated lattices]%
{Comaximal filter graphs in residuated lattices}
\author[Surdive Atamewoue* \and Hugue Tchantcho** ]%
{Surdive Atamewoue* \and Hugue Tchantcho*}

\newcommand{\acr}{\newline\indent}

\address{\llap{*\,} Department of Mathematics\acr
                  Advanced Teacher Training College \acr
                  University of Yaounde 1\acr
                   Yaounde\acr
                   Cameroon}
\email{surdive.atamewoue@univ-yaounde1.cm}

\address{\llap{*\,}Department of Mathematics\acr
	Faculty of Science \acr
	University of Yaounde 1\acr
	Yaounde\acr
	Cameroon}
\email{thuguesfr@yahoo.fr}


\thanks{The first author gratefully acknowledges the AIMS Cameroon Research Centre for the research stay they generously offered, which provided the opportunity to initiate this work.}


\begin{abstract}
Consider $\mathcal{A}$ to be a commutative, integral and non-degenerate residuated lattice. In this work, we introduce the graph of comaximal filters on the residuated lattice $\mathcal{A}$. We will denote by $\mathcal{C}_f(\mathcal{A})$ this graph for which the set of vertices are proper filters of $\mathcal{A}$ which are not contained in the radical of $\mathcal{A}$, and the adjacency relation on vertices is given as: consider two filters $F$ and $G$, there are adjacent if and only if $F\veebar G$; the filter generated by $F\cup G$ is equal to $\mathcal{A}$. The elementary properties of this graph are provided, we establish a link between this comaximal filter graphs and the zero-divisor graphs. Furthermore, we investigate the chromatic number, the clique number, the planarity and the perfection of the graph $\mathcal{C}_f(\mathcal{A})$. We also briefly describe all the comaximal filter graphs constructed on specific residuated lattices of small size.
\end{abstract}

\maketitle

\section{Introduction}
Nowadays it is essential to establish rational logical systems as the foundational logic for processing uncertain informations. For this reason, a variety of non-classical logic systems have been extensively proposed and studied. In fact, non-classical logic has become a formal and valuable tool in computer science for handling uncertain and fuzzy informations. On the other hand, several logical algebras have been proposed as the semantic systems for non-classical logic systems, such as residuated lattices, MV-algebras, BL-algebras, G\"{o}del algebras, MTL-algebras, Heyting algebras, and Boolean-algebras \cite{EG01, GJKO07, H98, WD39, T01}, among others. Among these logical algebras, residuated lattices are fundamental and important algebraic structures because all these other logical algebras are specific cases of residuated lattices. The filter theory of logical algebras plays a significant role in the study of these algebras and the completeness of the corresponding non-classical logics. All of this justifies why the notion of a filter is extensively studied on residuated lattices. However, there are several approaches to this notion: classical filters, implicative filters, fantastic filters, obstinate filters, positive implicative filters, pure filters, easy filters (the reader is referred to \cite{BP15} for more on filters). But in this work, we will use the classical definition of filter.

Many authors have studied graph theory in relation to algebraic structures. A graph is a set of vertices connected by edges, used to model relationships between objects. The connection between graphs and algebraic structures is not a recent development. Over the past thirty years, using graphs to explore algebraic structures has emerged as an exciting area of research. Beck (1988) introduced the concept of the zero-divisor graph of a commutative ring, although his work primarily focused on ring colorings \cite{B88}. After that, Anderson and Livingston in \cite{AL99} altered the notion and his definition does not consider zero of the ring as a vertex of the graph. More recently the zero-divisor graph has been extended to other algebraic structures; as on semirings with the works of DeMeyer et al. \cite{DMS02}, DeMeyer and DeMeyer \cite{DD05}, and Miller et al. \cite{MTT15}; as on residuated lattices where Torkzadeh et al. \cite{TAB14} introduced the notion of zero-divisors of a
nonempty subset of a residuated lattice and he used that set to define the residuated graph.  In \cite{SB95}, Sharma and Bhatwadekar, introduced the co-maximal graph on a commutative ring $R$, where the vertices are the elements of $R$, and two distinct vertices $x$ and $y$ are adjacent if and only if $Rx + Ry = R$. Maimani et al. \cite{MSSY08} have also investigated this class of graphs, with a particular emphasis on the subgraph induced by the non-unit elements of the ring. In 2011, an approach to the study of the interconnection between graphs theory and commutative rings was proposed in the works of Behboodi and Rakeei \cite{BR11}, where they  defined a new graph structure on commutative
rings that used ideals instead of elements of a ring. Following in the same direction as Behboodi et al., Ye and Wu \cite{YW12} introduced in 2012 the notion of the co-maximal ideal graph of a commutative ring, where the vertices are the ideals of the ring rather than the elements of the ring, as in Sharma and Bhatwadekar's approach. 

In this work, we focus on the connection between graph theory and the theory of residuated lattices. Our goal is to introduce and study the comaximal filter graphs constructed from residuated lattices. The elementary properties of this graph are provided, we establish a link between this comaximal filter graphs and the zero-divisor graphs. Furthermore, we investigate the chromatic number, the clique number, the planarity and the perfection of these  comaximal filter graphs. We also briefly describe all the comaximal filter graphs constructed on specific residuated lattices of size less than or equal to 10. 

\section{Preliminaries}
\subsection{On residuated lattices}\leavevmode
\begin{defn}\cite{WD39}
A commutative bounded integral residuated lattice is an algebra $\mathcal{A} = (A; \wedge,\vee, \odot, \rightarrow, 0, 1)$ of type (2, 2, 2, 2, 0, 0) satisfying the following conditions:
\begin{description}
	\item[i] $(A; \odot, 1)$ is a commutative monoid;
	\item[ii] $(A; \wedge, \vee, 0, 1)$ is a bounded lattice;
	\item[iii] $x \odot y \leq z$ if and only if  $x\leq y\rightarrow z$, for any $x,y,z\in A$.
\end{description}
\end{defn}
In what follows, by a residuated lattice we will mean a commutative bounded integral residuated lattice.
We define the unary operation (negation) ``$\neg$" on $\A$ by $\neg x := x \rightarrow 0$ for any $x \in A$.\\
A residuated lattice $\A$ is called:

an \textit{MTL-algebra} \cite{EG01} if $\A$ satisfies the identity of pre-linearity, that is 
\begin{description}
	\item[iv] $(x \rightarrow y) \vee (y \rightarrow x) = 1$ (pre-linearity);
\end{description}

\textit{involutive} if $\A$ satisfies the identity of double negation, that is
\begin{description}
	\item[v] $\neg\neg x=x$ (double negation);
\end{description}

\textit{divisible} if $\A$ satisfies the identity of divisibility, that is
\begin{description}
	\item[vi] $x\wedge y= x\odot (x \rightarrow y)$ (divisibility);
\end{description}

a \textit{BL-algebra} \cite{H98} if $\A$ satisfies both (\textbf{iv}) and (\textbf{vi});

a \textit{G\"{o}del algebra} \cite{ZX10} if $\A$ is a \textit{BL-algebra} satisfying the identity of idempotence, that is 
\begin{description}
	\item[vii] $x\odot x=x$ (idempotence);
\end{description}

an \textit{MV-algebra} \cite{C58, M07} if $\A$ is an involutive BL-algebra;

a \textit{Heyting algebra} \cite{VDCK10} if the operations ``$\odot$" and ``$\wedge$" coincide on $\A$, that is, $\A$ satisfies (\textbf{vii}).

\begin{defn}\cite{SZ06}
Let $\A$ be a residuated lattice and  $\emptyset \neq F \subseteq A$, then $F$ is called a \textit{filter} of $\A$
if for any $x, y \in F$ and $z \in A$:
\begin{description}
	\item[i] $x\odot y \in F$;
	\item[ii] $x \leq z$ implies $z \in F$. (upward closed)
\end{description}
\end{defn}
\begin{pro}\cite{T01} 
Let $\A$ be a residuated lattice.
If $\emptyset \neq F \subseteq A$ then, $F$ is a filter of $\A$ if and only if for any $x, y \in A$;
$x \in F$ and $x\rightarrow y \in F$ imply $y \in F$. That means $F$ is a \textit{deductive system} of $\A$.
\end{pro}
We denote by $\mathcal{F}(\A)$ the set of all filters of a residuated lattice $\A$. It is well known that $(\mathcal{F}(\A), \subseteq)$ is a complete lattice in which infima are equal to the set intersections \cite{GJKO07}.

Let $X\subseteq  A$, we denote by $\langle X\rangle$ the filter of $\A$ generated by $X$. We have
\[\langle X \rangle = \{a\in A: a\geq x_1\odot x_2\odot...\odot x_n, ~ \textrm{where}~ n \in \mathbb{N},~~ x_1,..., x_n \in X\}. \]
Let $X$ and $Y$ be two subsets of $\A$, we denote by $X\veebar Y$ the filter of $\A$ generated by $X\cup Y$. Then 
\[X\veebar Y = \{a\in A: a\geq x\odot y, ~ \textrm{where}~ x \in X,  y \in Y\}. \]
A filter $F$ of $\A$ is called maximal when for any filter  $G$ of $\A$, if $F \subseteq G$, then $F = G$ or $G = A$. We denote by $Max(\A)$ the set of all proper maximal filter of $\A$. A residuated lattice that has only one maximal filter is called local. 
\begin{pro}\label{prop0}\cite{JG09}
Let $\A$ be a finite residuated lattice. Any filter of $\A$ is generated by an idempotent less than or equal to 1.
\end{pro}

\begin{defn}\cite{KO01}
Let $\A$ be a residuated lattice. The intersection of all maximal filters of $\A$ is called the radical of $\A$ and will be denoted by $Rad(\A)$. It is obvious that $Rad(\A) \in \mathcal{F}(\A)$.
\end{defn}
\subsection{On graphs theory}\leavevmode

Let's now review some facts and notations on graphs theory relevant to this paper. All the graph-related concepts mentioned in this subsection are taken from \cite{B79, BM08}, which also covers graphs in much greater detail. We therefore refer the interested reader to these sources for further learning on the subject.

A graph \( \Gamma \) is defined as an ordered pair \( (V(\Gamma), E(\Gamma)) \), where \( V(\Gamma) \) denotes the set of vertices and \( E(\Gamma) \) the set of edges, such that \( E(\Gamma) \cap V(\Gamma) = \emptyset \). A graph is said to be \textit{connected} if, for every pair of distinct vertices, there exists a path linking them. Let \( u \) and \( v \) be two vertices of $\Gamma$, the length of the shortest path connecting them is called the \textit{distance} between \( u \) and \( v \), and denoted by \( d(u, v) \). If no such path exists, we set \( d(u, v) = \infty \). The supremum of the distance between all vertices of $\Gamma$ is called the \textit{diameter} of $\Gamma$, and denoted by $diam(\Gamma)$.
 
Assume that $\Gamma$ contains cycles, then the $girth$ of $\Gamma$ (symbolized by $girth(\Gamma)$) is the length of its shortest cycle. A graph with no cycles has infinite $girth$. $\Gamma$ is said to be \textit{$s$-partite} if its vertex set can be partitioned into \( s \) disjoint subsets such that no edge has both of its endpoints within the same subset.  

A \textit{complete $s$-partite} is an \( s \)-partite graph in which any two vertices belonging to different subsets are adjacent. The complete bipartite graph with part sizes $p$ and $q$ is denoted by $\mathcal{K}_{p,q}$. A \textit{complete graph} is one in which every two distinct vertices are adjacent, meaning that an edge connects each possible pair of vertices. We use $\mathcal{K}_n$ to
denote the complete graph with $n$ vertices. \textit{Path graph} on $n$ vertices, is denoted by $\mathcal{P}_n$.

When a subgraph induced on a subset $X \subseteq V (\Gamma)$ is complete, then $X$ is called a \textit{clique}. The \textit{clique number} of a graph \( \Gamma \), denoted by \( \omega(\Gamma) \), is defined as the number of vertices in a largest complete subgraph (clique) of \( \Gamma \).

A \textit{coloring} of a graph \( \Gamma \) is a mapping that assigns a color (from a given set) to each vertex of \( \Gamma \), in such a way that no two adjacent vertices share the same color. When exactly \( c \) colors are used in such an assignment, it is called an \( c \)-coloring. A graph \( \Gamma \) is said to be \( c \)-colorable if there exists an \( c \)-coloring of \( \Gamma \). The smallest such \( c \) for which \( \Gamma \) admits an \( c \)-coloring is called the chromatic number of \( \Gamma \), denoted by \( \chi(\Gamma) \). 

Two graphs are said to be \textit{isomorphic} if there exists a bijection between their vertex sets that preserves adjacency.  Determining whether two graphs are isomorphic is a well-known computational problem with significant applications in chemistry, network analysis, and computer science.

A graph is said to be \textit{planar} if it can be drawn in the plane without any of its edges crossing, except at their endpoints. In other words, a planar graph admits an embedding in the plane such that its edges intersect only at vertices. Classic results such as Kuratowski's and Wagner's theorems \cite{K30} provide powerful characterizations of planarity. 

\section{The comaximal filter graphs in residuated lattices}\label{Sect2}

Let us now focus on the comaximal filter graphs in residuated lattices. Let $\A=(A; \wedge, \vee, \odot, \rightarrow, 0, 1)$ be a residuated lattice. The comaximal filter graph on $\A$, denoted by  $\mathcal{C}_f(\mathcal{A})$, is defined as the simple graph with:\\
\indent - \textbf{Vertex set:} $V(\mathcal{C}_f(\mathcal{A}))$ are proper filters of $\mathcal{A}$ which are not contained in the radical of $\mathcal{A}$.\\  
\indent - \textbf{Edge set:} \(E(\mathcal{C}_f(\mathcal{A})) = \{\{F, G\} \mid F \ne G,\, F, G \in \mathcal{F}(\A),\, \text{and } F \veebar G = A\}\).

\begin{pro}\label{prop1} Let $\mathcal{A}=(A,\wedge,\vee, \odot, \rightarrow, 0, 1)$ be a residuated lattice. Then the following assertions hold.
	\begin{enumerate}
		\item $\mathcal{C}_f(\mathcal{A})$ is the empty graph if and only if $\A$ is a local residuated lattice.
		\item If $|Max(\A)| \geq 2$. Then for any proper filter $F$ such that $F\nsubseteq Rad(\A)$, $F$ is necessarily a vertex of $\CC_f(\A)$.
	\end{enumerate}
\end{pro}
\begin{proof}
(1) Clear since any proper filter should be contains in the only maximal filter \cite{C06}.\\
(2) $F\nsubseteq Rad(\A)$ provide the existence of a maximal filter $M$ not containing $F$. The maximal property of $M$ combine with $M\subsetneq M\veebar F$, lead to $M\veebar F = A$. So $\{F, M\}$ is an edge of $\CC_f(\A)$.
\end{proof} 
The Proposition \ref{prop1} highlights the fact that the study of graphs $\CC_f(\A)$ where $\A$ is local cannot be interesting. For this reason we focus on graphs construct on non-local residuated lattices.
\begin{ex}\label{ex1}
Let $A = \{0, a, b, c, d, e, f, 1\}$ be a set and $\A = (A; \wedge, \vee, \odot, \rightarrow, 0, 1)$ be a residuated lattice in which;\\ the commutative operation $``\odot"$ is given by\quad
$\begin{tabular}{c|c|c|c|c|c|c|c|c}
	$\odot$ & 0 & a & b & c & d & e & f & 1 \\ 	\hline 
	0 & 0 & 0 & 0 & 0 & 0 & 0 & 0 & 0 \\  	\hline 
	a & 0 & a & 0 & a & 0 & a & 0 & 1 \\  	\hline
	b & 0 & 0 & b & b & 0 & 0 & b & b \\  	\hline
	c & 0 & a & b & c & 0 & a & b & c \\  	\hline 
	d & 0 & 0 & 0 & 0 & d & d & d & d \\  	\hline 
	e & 0 & a & 0 & a & d & e & d & e \\  	\hline
	f & 0 & 0 & b & b & d & d & f & f \\    \hline
    1 &	0 & a & b & c & d & e & f & 1            
\end{tabular}$,\\\\
the order in $\A$, that is $``\leq"$, is given by the matrix that contains the value 1 in row $x$ and column $y$ if $x\leq y$, and 0 otherwise. \quad\quad\quad\quad\quad 
\begin{tabular}{c|c|c|c|c|c|c|c|c}
	$\leq$ & 0 & a & b & c & d & e & f & 1 \\ 	\hline 
	0 & 1 & 1 & 1 & 1 & 1 & 1 & 1 & 1 \\  	\hline 
	a & 0 & 1 & 0 & 1 & 0 & 1 & 0 & 1 \\  	\hline
	b & 0 & 0 & 1 & 1 & 0 & 0 & 1 & 1 \\  	\hline
	c & 0 & 0 & 0 & 1 & 0 & 0 & 0 & 1 \\  	\hline 
	d & 0 & 0 & 0 & 0 & 1 & 1 & 1 & 1 \\  	\hline 
	e & 0 & 0 & 0 & 0 & 0 & 1 & 0 & 1 \\    \hline
	f & 0 & 0 & 0 & 0 & 0 & 0 & 1 & 1 \\    \hline
	1 & 0 & 0 & 0 & 0 & 0 & 0 & 0 & 1
\end{tabular}, \\\\
the operation $``\rightarrow"$ is given by $x \rightarrow y = \vee\{a \in A~|~x\odot a \leq y\}$ for any $x, y \in A$. Through simple verification calculations, it is shown that $\A$ is an MV-algebra.\\ 
The set of filters of $\A$ is given by $\mathcal{F}(\A)=\{F_1=\{1\}, F_2=\{c,1\}, F_3=\{e,1\}, F_4=\{f,1\}, F_5=\{a, c, e, 1\}, F_6=\{b, c, f, 1\}, F_7=\{d, e, f, 1\},  F_8=A\}$.
Since $Rad(\A)=\{1\}=F_1$, then for the graph $\CC_f(\A)$, the set of vertices in given by $V(\CC_f(\A))=\{F_2, F_3, F_4, F_5, F_6, F_7\}$ and we then represent the graph as follows:\\
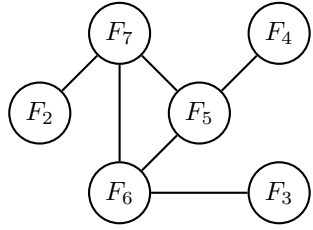
\begin{center}
	
\begin{minipage}{4cm}
	\begin{tikzpicture}[node distance={15mm}, thick, main/.style = {draw, circle}] 
	\node[main] (1) {$F_2$}; 
	\node[main] (2) [above right of=1] {$F_7$}; 
	\node[main] (3) [below right of=1] {$F_6$}; 
	\node[main] (4) [above right of=3] {$F_5$}; 
	\node[main] (5) [above right of=4] {$F_4$}; 
	\node[main] (6) [below right of=4] {$F_3$}; 
	\draw (1) -- (2); 
	\draw (2) -- (3); 
	\draw (2) -- (4); 
	\draw (3) -- (4); 
	\draw (3) -- (6); 
	\draw (4) -- (5); 
	\end{tikzpicture} 
\end{minipage}
\end{center}	
\end{ex}

\begin{pro}
Let $\CC_f(\A)$ be a comaximal filter graph on $\A$. Then $|V(\CC_f(\A))|\leq |\mathcal{F}(\A)|-2$ and this bound is attainable.	
\end{pro}
\begin{proof}
The definition of the graph and the Example \ref{ex1}, are sufficient to proof this.	
\end{proof}

In a residuated lattice, a comaximal filter is a subset of the lattice that satisfies certain conditions related to the lattice operations. There are well-known characterizations of comaximal filters \cite{RD22} in residuated lattices that closely relate this notion to that of zero-divisors. 
\begin{pro}\cite{RD22, RD25}\label{Prop1}
Let $\A$ be a residuated lattice and $F, G$ two proper filters of $\A$. The following assertions are equivalent:
\begin{enumerate}
	\item $F$ and $G$ are comaximal;
	\item there exist $x \in F$ and $y \in G$ such that $x \odot y = 0$.
\end{enumerate}
\end{pro}

\begin{rem}\label{rem1}
	Two different maximal filters of a residuated lattice are comaximal.
\end{rem}

In \cite{GY20}, Gan and Yang introduced the notion of zero-divisor graph of an MV-algebra. This concept can be extended to any residuated lattice. Let $\A=(A; \wedge, \vee, \odot, \rightarrow, 0, 1)$, and \( (A, \odot, 0) \) its associated commutative semigroup. The \textit{zero-divisor graph} of \( \A \), denoted \( \Gamma(\A) \), is defined as the simple graph with vertex set \( V(\Gamma(\A)) = \{x \in A \mid \exists y \in A \setminus \{0\} \text{ such that } x \odot y = 0\} \), and edge set \( E(\Gamma(\A)) = \{\{x, y\} \mid x \neq y,\ x, y \in A,\ \text{and } x \odot y = 0\} \).

The upcoming results and observations establish a connection between zero-divisor graphs and comaximal filter graphs in the context of residuated lattices.

\begin{rem}From Proposition \ref{Prop1}, it follows that  for a residuated lattice, the comaximal filter graph is a subgraph of the zero-divisor graph in the sense of Gan and Yang \cite{GY20}. 
\end{rem}

Let $\A=(A,\wedge,\vee,\odot,\rightarrow,0,1)$ be a residuated lattice. We define by $Z(\A)\setminus\{0\}$ the set of non zero-divisors of $\A$. That is, an element $x\in A$ for which there exists $y\in A\setminus\{0\}$ such that $x\odot y=0$.

\begin{thm}
Let $\A=(A,\wedge,\vee,\odot,\rightarrow,0,1)$ be a finite residuated lattice. If $|Z(\A)\setminus\{0\}|=|V(\CC_f(\A))|$, then the comaximal filter graph $\CC_f(\A)$ and the zero-divisor graph $\Gamma(\A)$ are isomorphic. 
\end{thm}
\begin{proof} The map $\varphi$ defined as follows
	\[\begin{array}{rcl}
		\varphi: Z(\A)\setminus\{0\} &\to& V(\CC_f(\A))\\
		x &\mapsto & \varphi(x)=\langle x\rangle
	\end{array}\] is the one for the graphs isomorphism.
	
	In fact, considering the fact that composing two different elements of a filter using operation ``$\odot$" cannot result in $0$, thus two distinct elements of $|Z(\A)\setminus\{0\}|$ cannot generate the same filter. So, if $\{x,y\}\in E(\Gamma(\A))$, then by Proposition \ref{Prop1} $\{\langle x\rangle, \langle y\rangle\}\in E(\CC_f(\A))$. Conversely, let $\{F, G\}\in E(\CC_f(\A))$. Since every filter in a finite residuated lattice is generated by an idempotent less than or equal to 1 \ref{prop0}, then assuming $F = \langle x\rangle$ and $G = \langle y\rangle$, we easily obtain that $x\odot y = 0$. Therefore, $\{x, y\}\in E(\Gamma(\A))$.	
\end{proof}

\begin{pro}
Let $\A$ be a finite residuated lattice. The comaximal filter graph $\CC_f(\A)$ is isomorphic to a subgraph of the zero-divisor graph $\Gamma(\A)$.
\end{pro}
\begin{proof}
For any finite residuated lattice, using Proposition \ref{Prop1} we conclude that   
$|V(\mathcal{C}_f(\mathcal{A}))| \leq |Z(\mathcal{A}) \setminus \{0\}|$. 
\end{proof}

Thus comaximal filter graphs are not zero-divisor graphs. This is why it is appropriate to analyze and study this graph to determine which properties of zero-divisor graphs are preserved.

\begin{ex}\label{ex2}
	Let $\A$ be a residuated lattice with operations table:\\
	\indent \indent  $\begin{tabular}{c|c|c|c|c|c|c|c|c|c|c}
	$\leq$ & 0 & a & b & c & d & e & f & g & h & 1 \\ 	\hline 
	0 & 1 & 1 & 1 & 1 & 1 & 1 & 1 & 1 & 1 & 1 \\  	\hline 
	a & 0 & 1 & 0 & 1 & 1 & 1 & 1 & 1 & 1 & 1 \\  	\hline 
	b & 0 & 0 & 1 & 1 & 0 & 1 & 0 & 1 & 0 & 1 \\  	\hline
	c & 0 & 0 & 0 & 1 & 0 & 1 & 0 & 1 & 0 & 1 \\  	\hline
	d & 0 & 0 & 0 & 0 & 1 & 1 & 1 & 1 & 1 & 1 \\  	\hline 
	e & 0 & 0 & 0 & 0 & 0 & 1 & 0 & 1 & 0 & 1 \\  	\hline 
	f & 0 & 0 & 0 & 0 & 0 & 0 & 1 & 1 & 1 & 1 \\    \hline
	g & 0 & 0 & 0 & 0 & 0 & 0 & 0 & 1 & 0 & 1 \\    \hline
	h & 0 & 0 & 0 & 0 & 0 & 0 & 0 & 0 & 1 & 1 \\    \hline
	1 & 0 & 0 & 0 & 0 & 0 & 0 & 0 & 0 & 0 & 1    
	\end{tabular}$\quad; \quad   $\begin{tabular}{c|c|c|c|c|c|c|c|c|c|c}
	$\odot$ & 0 & a & b & c & d & e & f & g & h & 1 \\ 	\hline 
	0 & 0 & 0 & 0 & 0 & 0 & 0 & 0 & 0 & 0 & 0 \\  	\hline 
	a & 0 & 0 & 0 & 0 & 0 & 0 & 0 & 0 & a & a \\  	\hline
	b & 0 & 0 & b & b & 0 & b & 0 & b & 0 & b \\  	\hline
	c & 0 & 0 & b & b & 0 & b & 0 & b & a & c \\  	\hline 
	d & 0 & 0 & 0 & 0 & 0 & 0 & a & a & d & d \\  	\hline 
	e & 0 & 0 & b & b & 0 & b & a & c & d & e \\  	\hline
	f & 0 & 0 & 0 & 0 & a & a & d & d & f & f \\    \hline
	g & 0 & 0 & b & b & a & c & d & e & f & g \\    \hline
	h & 0 & a & 0 & a & d & d & f & f & h & h \\    \hline
	1 &	0 & a & b & c & d & e & f & g & h & 1            
	\end{tabular}$\\\\
The operation $``\rightarrow"$ is defined as in Example \ref{ex1}.\\ 
We obtain that $V(\CC_f(\A))=\{F_1=\{h,1\}=\langle h\rangle, F_2=\{b,c,e,g,1\}=\langle b\rangle\}$, and $V(Z(\A)\setminus\{0\})=\{a,b,c,d,e,f,g,h\}$.\\
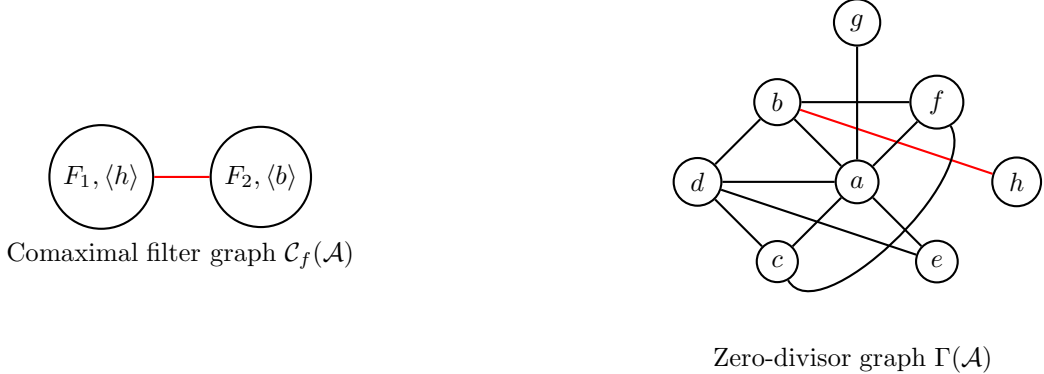
\begin{minipage}{6cm}
	\centering
\begin{tikzpicture}[node distance={15mm}, thick, main/.style = {draw, circle}] 
	\node[main] (3) [below right of=1] {$F_1,\langle h\rangle$}; 
	\node[main] (6) [below right of=4] {$F_2,\langle b\rangle$}; 
	\draw[red] (3) -- (6); 
	\end{tikzpicture} \\
	Comaximal filter graph $\CC_f(\A)$
\end{minipage}\hfill
\begin{minipage}{6cm}
	\centering
	\begin{tikzpicture}[node distance={15mm}, thick, main/.style = {draw, circle}] 
	\node[main] (1) {$d$}; 
	\node[main] (2) [above right of=1] {$b$}; 
	\node[main] (3) [below right of=1] {$c$}; 
	\node[main] (4) [above right of=2] {$g$}; 
	\node[main] (5) [below right of=2] {$a$}; 
	\node[main] (6) [below right of=5] {$e$}; 
	\node[main] (7) [above right of=5] {$f$}; 
	\node[main] (8) [below right of=7] {$h$}; 
	\draw (1) -- (2); 
	\draw (1) -- (3); 
	\draw (1) -- (5); 
	\draw (1) -- (6); 
	\draw (2) -- (5); 
	\draw (2) -- (7);
	\draw[red] (2) -- (8); 
	\draw (3) to [out=300, in=300, looseness=0.7] (7); 
	\draw (4) -- (5); 
	\draw (5) -- (6); 
	\draw (5) -- (7); 
	\draw (3) -- (5); 
	\end{tikzpicture}\\ 
	Zero-divisor graph $\Gamma(\A)$
	\bigskip
\end{minipage}\\
\noindent By scrutinizing the graphs above, it is clear that $\CC_f(\A)$ is a subgraph of $\Gamma(\A)$, yet it appears to share no properties with graph $\Gamma(\A)$.
\end{ex}

\begin{pro}
	Isomporhic residuated lattices induce comaximal filter graphs isomorphism. 
\end{pro}
\begin{proof}
	Two residuated lattices are said to be isomorphic if there exists a bijective map between them that preserves the lattice operations as well as the monoidal operation and the residual. So, their filters should be isomorphic. 
\end{proof}
The converse of this proposition is not true in general. See the below example.
\begin{ex}
	Let $\mathcal{G}_1$ and $\mathcal{G}_2$ be two G\"{o}del-algebras, with the operations tables given by:\\\\
	$\mathcal{G}_1$: $\begin{tabular}{c|c|c|c|c|c|c|c|c|c|c}
	$\leq$ & 0 & a & b & c & d & e & f & g & h & 1 \\ 	\hline 
	0 & 1 & 1 & 1 & 1 & 1 & 1 & 1 & 1 & 1 & 1 \\  	\hline 
	a & 0 & 1 & 0 & 1 & 1 & 1 & 1 & 1 & 1 & 1 \\  	\hline 
	b & 0 & 0 & 1 & 1 & 0 & 1 & 0 & 0 & 1 & 1 \\  	\hline
	c & 0 & 0 & 0 & 1 & 0 & 1 & 0 & 0 & 1 & 1 \\  	\hline
	d & 0 & 0 & 0 & 0 & 1 & 1 & 0 & 1 & 0 & 1 \\  	\hline 
	e & 0 & 0 & 0 & 0 & 0 & 1 & 0 & 0 & 0 & 1 \\  	\hline 
	f & 0 & 0 & 0 & 0 & 0 & 0 & 1 & 1 & 1 & 1 \\    \hline
	g & 0 & 0 & 0 & 0 & 0 & 0 & 0 & 1 & 0 & 1 \\    \hline
	h & 0 & 0 & 0 & 0 & 0 & 0 & 0 & 0 & 1 & 1 \\    \hline
	1 & 0 & 0 & 0 & 0 & 0 & 0 & 0 & 0 & 0 & 1    
	\end{tabular}$;\quad and \quad	
	 $\begin{tabular}{c|c|c|c|c|c|c|c|c|c|c}
	$\odot$ & 0 & a & b & c & d & e & f & g & h & 1 \\ 	\hline 
	0 & 0 & 0 & 0 & 0 & 0 & 0 & 0 & 0 & 0 & 0 \\  	\hline 
	a & 0 & a & 0 & a & a & a & a & a & a & a \\  	\hline
	b & 0 & 0 & b & b & 0 & b & 0 & 0 & b & b \\  	\hline
	c & 0 & a & b & c & a & c & a & a & c & c \\  	\hline 
	d & 0 & a & 0 & a & d & d & a & d & a & d \\  	\hline 
	e & 0 & a & b & c & d & e & a & d & c & e \\  	\hline
	f & 0 & a & 0 & a & a & a & f & f & f & f \\    \hline
	g & 0 & a & 0 & a & d & d & f & g & f & g \\    \hline
	h & 0 & a & b & c & a & c & f & f & h & h \\    \hline
	1 &	0 & a & b & c & d & e & f & g & h & 1            
	\end{tabular}$ \bigskip \\
    $\mathcal{G}_2$: $\begin{tabular}{c|c|c|c|c|c|c|c|c|c|c}
	$\leq$ & 0 & a & b & c & d & e & f & g & h & 1 \\ 	\hline 
	0 & 1 & 1 & 1 & 1 & 1 & 1 & 1 & 1 & 1 & 1 \\  	\hline 
	a & 0 & 1 & 0 & 1 & 1 & 1 & 1 & 1 & 1 & 1 \\  	\hline 
	b & 0 & 0 & 1 & 1 & 0 & 1 & 0 & 1 & 0 & 1 \\  	\hline
	c & 0 & 0 & 0 & 1 & 0 & 1 & 0 & 1 & 0 & 1 \\  	\hline
	d & 0 & 0 & 0 & 0 & 1 & 1 & 1 & 1 & 1 & 1 \\  	\hline 
	e & 0 & 0 & 0 & 0 & 0 & 1 & 0 & 1 & 0 & 1 \\  	\hline 
	f & 0 & 0 & 0 & 0 & 0 & 0 & 1 & 1 & 1 & 1 \\    \hline
	g & 0 & 0 & 0 & 0 & 0 & 0 & 0 & 1 & 0 & 1 \\    \hline
	h & 0 & 0 & 0 & 0 & 0 & 0 & 0 & 0 & 1 & 1 \\    \hline
	1 & 0 & 0 & 0 & 0 & 0 & 0 & 0 & 0 & 0 & 1    
	\end{tabular}$;\quad and \quad	
	$\begin{tabular}{c|c|c|c|c|c|c|c|c|c|c}
	$\odot$ & 0 & a & b & c & d & e & f & g & h & 1 \\ 	\hline 
	0 & 0 & 0 & 0 & 0 & 0 & 0 & 0 & 0 & 0 & 0 \\  	\hline 
	a & 0 & a & 0 & a & a & a & a & a & a & a \\  	\hline
	b & 0 & 0 & b & b & 0 & b & 0 & b & 0 & b \\  	\hline
	c & 0 & a & b & c & a & c & a & c & a & c \\  	\hline 
	d & 0 & a & 0 & a & d & d & d & d & d & d \\  	\hline 
	e & 0 & a & b & c & d & e & d & e & d & e \\  	\hline
	f & 0 & a & 0 & a & d & d & f & f & f & f \\    \hline
	g & 0 & a & b & c & d & e & f & g & f & g \\    \hline
	h & 0 & a & 0 & a & d & d & f & f & h & h \\    \hline
	1 &	0 & a & b & c & d & e & f & g & h & 1            
	\end{tabular}$ \\\\		
The operations $``\rightarrow"$ is defined as in Example \ref{ex1}.\\\\
The comaximal filter graphs are quickly obtain as:\\ 	
$V(\CC_f(\mathcal{G}_1))=\{F_1=\{g,1\}, F_2=\{d,e,g,1\}, F_3=\{f, g, h, 1\}, F_4=\{b,c,e,h,1\}, F_5=\{a, c, d, e, f, g, h, 1\} \}$	and $E(\CC_f(\mathcal{G}_1))=\{F_1-F_4, F_2-F_4, F_3-F_4, F_4-F_5\}$ for $\mathcal{G}_1$.\\\\
And\\\\
$V(\CC_f(\mathcal{G}_2))=\{F_1'=\{h,1\}, F_2'=\{f,g,h,1\}, F_3'=\{b,c,e, g, 1\}, F_4'=\{d,e,f, g, h,1\}, F_5'=\{a, c, d, e, f, g, h, 1\} \}$	and $E(\CC_f(\mathcal{G}_2))=\{F_1'-F_3', F_2'-F_3', F_3'-F_4', F_3'-F_5'\}$ for $\mathcal{G}_2$.\\\\
The map 	\[\begin{array}{rcl}
\varphi:V(\CC_f(\mathcal{G}_1)) &\to& V(\CC_f(\mathcal{G}_2))\\
F_1 &\mapsto & \varphi(F_1)=F_1'\\
F_2 &\mapsto & \varphi(F_2)=F_2'\\
F_3 &\mapsto & \varphi(F_3)=F_4'\\
F_4 &\mapsto & \varphi(F_4)=F_3'\\
F_5 &\mapsto & \varphi(F_5)=F_5'
\end{array}\] define a graph isomorphism between $\CC_f(\mathcal{G}_1)$ and $\CC_f(\mathcal{G}_2)$.\\ 
\begin{minipage}{4cm}
	\centering
	\begin{tikzpicture}[node distance={15mm}, thick, main/.style = {draw, circle}] 
	\node[main] (2) [above right of=1] {$F_1$}; 
	\node[main] (3) [below right of=1] {$F_5$}; 
	\node[main] (4) [above right of=3] {$F_4$}; 
	\node[main] (5) [above right of=4] {$F_2$}; 
	\node[main] (6) [below right of=4] {$F_3$}; 
	\draw (4) -- (2); 
	\draw (4) -- (3); 
	\draw (4) -- (5); 
	\draw[red] (4) -- (6); 
	\end{tikzpicture} \\
	Comaximal filter graph $\CC_f(\mathcal{G}_1)$
\end{minipage}\hfill
\begin{minipage}{4cm}
	\centering
	\begin{tikzpicture}[node distance={15mm}, thick, main/.style = {draw, circle}] 
	\node[main] (2) [above right of=1] {$F_1'$}; 
	\node[main] (3) [below right of=1] {$F_5'$}; 
	\node[main] (4) [above right of=3] {$F_3'$}; 
	\node[main] (5) [above right of=4] {$F_2'$}; 
	\node[main] (6) [below right of=4] {$F_4'$}; 
	\draw (4) -- (2); 
	\draw (4) -- (3); 
	\draw (4) -- (5); 
	\draw[red] (4) -- (6); 
	\end{tikzpicture}\\ 
	Comaximal filter graph $\CC_f(\mathcal{G}_2)$\\
	\bigskip
\end{minipage}

\noindent By routine calculation, we find out that the two G\"{o}del algebra $\mathcal{G}_1$ and $\mathcal{G}_2$ are not isomporhic.
\end{ex}

\begin{rem}\label{Rem2}
	One can find residuated lattices of different sizes whose comaximal filter graphs constructed on each are isomorphic. Indeed, let $\A$ be a residuated lattice of size 9 with operations tables given by:\\
		\indent \indent  $\begin{tabular}{c|c|c|c|c|c|c|c|c|c}
	$\leq$ & 0 & a & b & c & d & e & f & g & 1 \\ 	\hline 
	0 & 1 & 1 & 1 & 1 & 1 & 1 & 1 & 1 & 1 \\  	\hline 
	a & 0 & 1 & 0 & 1 & 0 & 1 & 0 & 1 & 1 \\  	\hline 
	b & 0 & 0 & 1 & 1 & 0 & 0 & 1 & 1 & 1 \\  	\hline
	c & 0 & 0 & 0 & 1 & 0 & 0 & 0 & 1 & 1 \\  	\hline
	d & 0 & 0 & 0 & 0 & 1 & 1 & 1 & 1 & 1 \\  	\hline 
	e & 0 & 0 & 0 & 0 & 0 & 1 & 0 & 1 & 1 \\  	\hline 
	f & 0 & 0 & 0 & 0 & 0 & 0 & 1 & 1 & 1 \\    \hline
	g & 0 & 0 & 0 & 0 & 0 & 0 & 0 & 1 & 1 \\    \hline
	1 & 0 & 0 & 0 & 0 & 0 & 0 & 0 & 0 & 1    
	\end{tabular}$\quad; \quad   $\begin{tabular}{c|c|c|c|c|c|c|c|c|c}
	$\odot$ & 0 & a & b & c & d & e & f & g & 1 \\ 	\hline 
	0 & 0 & 0 & 0 & 0 & 0 & 0 & 0 & 0 & 0 \\  	\hline 
	a & 0 & a & 0 & a & 0 & a & 0 & a & a \\  	\hline
	b & 0 & 0 & b & b & 0 & 0 & b & b & b \\  	\hline
	c & 0 & a & b & c & 0 & a & b & c & c \\  	\hline 
	d & 0 & 0 & 0 & 0 & d & d & d & d & d \\  	\hline 
	e & 0 & a & 0 & a & d & e & d & e & e \\  	\hline
	f & 0 & 0 & b & b & d & d & f & f & f \\    \hline
	g & 0 & a & b & c & d & e & f & g & g \\    \hline
	1 &	0 & a & b & c & d & e & f & g & 1            
	\end{tabular}$\\\\
	The operation $``\rightarrow"$ is defined as in Example \ref{ex1}.\\
	We obtain that the graph $\CC_f(\A)$ define here by $V(\CC_f(\A))=\{F_1=\{c,g,1\},F_2=\{e,g,1\}, F_3=\{f,g,1\}, F_4=\{a,c,e,g,1\}, F_5=\{b,c,f,g,1\}, F_6=\{d,e,g,1\}\}$ and\\
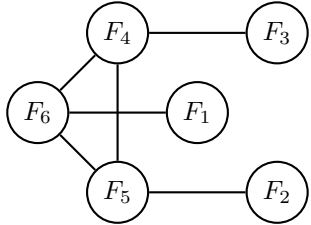
	\begin{center}
		
		\begin{minipage}{4cm}
			\begin{tikzpicture}[node distance={15mm}, thick, main/.style = {draw, circle}] 
			\node[main] (1) {$F_6$}; 
			\node[main] (2) [above right of=1] {$F_4$}; 
			\node[main] (3) [below right of=1] {$F_5$}; 
			\node[main] (4) [above right of=3] {$F_1$}; 
			\node[main] (5) [above right of=4] {$F_3$}; 
			\node[main] (6) [below right of=4] {$F_2$}; 
			\draw (1) -- (2); 
			\draw (2) -- (5); 
			\draw (2) -- (3); 
			\draw (1) -- (3); 
			\draw (1) -- (4); 
			\draw (3) -- (6); 
			\end{tikzpicture} 
		\end{minipage}
	\end{center}
is isomorphic to the comaximal filter graph obtain in Example \ref{ex1} on a residuated lattice of size 8.	
\end{rem}

We now give some particular properties of the comaximal filter graphs.
\begin{thm}
Let $\A$ be a residuated lattice. The comaximal filter graph $\CC_f(\A)$ is a simple, connected with diameter less than or equal to three.
\end{thm}
\begin{proof}
According to the definition of $\CC_f(\A)$, there are no parallel edges between two  comaximal filters. Since a filter cannot contain two zero-divisors, it follows that there are no loops. Therefore, \(\mathcal{C}_f(\mathcal{A})\) is a simple graph.
	
Let \( F \) and \( G \) be any two vertices in the graph \( \CC_f(\A) \). If the filter generate by $F\cup G$ satisfies \( F \veebar G = A \), then the distance between them is \( d(F, G) = 1 \). Now, consider the case where \( F \veebar G \ne A \). If there exists a filter \( H \) such that both \( F \veebar H = A \) and \( G \veebar H = A \), then \( d(F, G) = 2 \). In the absence of such a filter, we can identify two distinct maximal filter \( M_1 \) and \( M_2 \) satisfying \( F \veebar M_1 = A \) and \( G \veebar M_2 = A \) \cite{KO01}. Given that \( M_1 \ne M_2 \), it follows that \( M_1 \veebar M_2 = A \), which implies \( d(F, G) = 3 \). Since the diameter of a graph is defined as the greatest distance between any pair of vertices, we conclude that the diameter of \( \CC_f(\A) \) is at most 3. 

Moreover, this reasoning shows that the graph is connected which complete the proof.
\end{proof}

\begin{rem}
	The maximal diameter of comaximal filter graphs is realizable. Indeed, the diameter of the graph in Example \ref{ex1} is 3.
\end{rem}

\begin{thm}
Given a residuated lattice \( \A \) where $|Max(\A)|=2$ and \( \CC_f(\A) \), it follows that \( girth(\CC_f(\A)) \in \{4, \infty\} \).
\end{thm}
\begin{proof}
We may assume that \( Max(\A) = \{ M_1, M_2 \} \). Let \( F_1 \) and \( F_2 \) be two distinct filters of $\A$, neither of which is contained in \( Rad(\A) \). If \( F_1 \subseteq M_1 \) and \( F_2 \subseteq M_2 \), then it necessarily follows that: $ F_1 \veebar F_2 = A$. So $\CC_f(\A)$ is the complete graph \( \mathcal{K}_2 \), which is also denoted by \( \mathcal{K}_{1,1} \), or a star graph of the form \( \mathcal{K}_{1,n} \), where \( n \in \mathbb{N} \) and \( 2 \leq n \leq \infty \), or a complete bipartite graph \( \mathcal{K}_{m,n} \), where \( m, n \in \mathbb{N} \), with \( 2 \leq m \leq \infty \) and \( 2 \leq n \leq \infty \).\\
Thus, the proof is complete.
\end{proof}

\begin{ex}\label{ex3}
		Let $\A$ be a residuated lattice with operations table:\\
	\indent \indent  $\begin{tabular}{c|c|c|c|c|c|c|c|c|c}
	$\leq$ & 0 & a & b & c & d & e & f & g & 1 \\ 	\hline 
	0 & 1 & 1 & 1 & 1 & 1 & 1 & 1 & 1 & 1 \\  	\hline 
	a & 0 & 1 & 0 & 1 & 1 & 1 & 0 & 1 & 1 \\  	\hline 
	b & 0 & 0 & 1 & 1 & 0 & 1 & 1 & 1 & 1 \\  	\hline
	c & 0 & 0 & 0 & 1 & 0 & 1 & 0 & 1 & 1 \\  	\hline
	d & 0 & 0 & 0 & 0 & 1 & 1 & 0 & 0 & 1 \\  	\hline 
	e & 0 & 0 & 0 & 0 & 0 & 1 & 0 & 0 & 1 \\  	\hline 
	f & 0 & 0 & 0 & 0 & 0 & 0 & 1 & 1 & 1 \\    \hline
	g & 0 & 0 & 0 & 0 & 0 & 0 & 0 & 1 & 1 \\    \hline
	1 & 0 & 0 & 0 & 0 & 0 & 0 & 0 & 0 & 1    
	\end{tabular}$\quad; \quad   $\begin{tabular}{c|c|c|c|c|c|c|c|c|c}
	$\odot$ & 0 & a & b & c & d & e & f & g & 1 \\ 	\hline 
	0 & 0 & 0 & 0 & 0 & 0 & 0 & 0 & 0 & 0 \\  	\hline 
	a & 0 & a & 0 & a & a & a & 0 & a & a \\  	\hline
	b & 0 & 0 & b & b & 0 & b & b & b & b \\  	\hline
	c & 0 & a & b & c & a & c & b & c & c \\  	\hline 
	d & 0 & a & 0 & a & d & d & 0 & a & d \\  	\hline 
	e & 0 & a & b & c & d & e & b & c & e \\  	\hline
	f & 0 & 0 & b & b & 0 & b & f & f & f \\    \hline
	g & 0 & a & b & c & a & c & f & g & g \\    \hline
	1 &	0 & a & b & c & d & e & f & g & 1            
	\end{tabular}$\\\\
The operation $``\rightarrow"$ is defined as in Example \ref{ex1}.\\ 
We obtain that $V(\CC_f(\A))=\{F_1=\{d,e,1\}, F_2=\{f,g,1\}, F_3=\{a,c,d,e,g,1\}, F_4=\{b,c,e,f,g,1\}\}$ and the graph\\
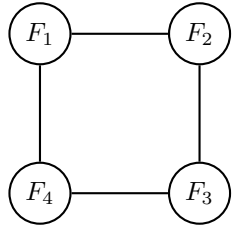
\begin{center}
	
	\begin{minipage}{4cm}
		\begin{tikzpicture}[node distance={15mm}, thick, main/.style = {draw, circle}] 
		\node[main] (2) [above right of=1] {$F_1$}; 
		\node[main] (3) [below right of=1] {$F_4$}; 
		\node[main] (5) [above right of=4] {$F_2$}; 
		\node[main] (6) [below right of=4] {$F_3$}; 
		\draw (2) -- (5); 
		\draw (2) -- (3); 
		\draw (6) -- (5); 
		\draw (3) -- (6); 
		\end{tikzpicture} 
	\end{minipage}
\end{center}
We observe that the $girth(\CC_f(\A))=4$.	
\end{ex}

\begin{rem}
	The $girth$ of the comaximal filters graphs in Example \ref{ex2} is $\infty$, since there graph is acyclic.
\end{rem}

As in the specific cases where \( Max(\A) = 1 \) or \( Max(\A) = 2 \), the $girth$ of the graph is well known, we will now focus on the $girth$ in the general case, that is, when \( Max(\A) \geq 3 \).

\begin{pro}
Let \( \A \) be a residuated lattice. If the number of maximal filters satisfies \( |Max(\A)| \geq 3 \), then \( \CC_f(\A) \) contains an induced subgraph isomorphic to the complete graph \( \mathcal{K}_3 \). Consequently, \( \CC_f(\A) \) must contain at least one cycle, and the $girth(\CC_f(\A)) = 3 $.
\end{pro}
\begin{proof}
Assume that \( |Max(\A)| \geq 3 \), meaning the residuated lattice \( \A \) possesses at least three distinct maximal filters, denoted \( F_1, F_2, \) and \( F_3 \). Let us examine the subgraph of \( \CC_f(\A) \) induced by the set \( \{F_1, F_2, F_3\} \). This subgraph forms a complete graph (since Remark \ref{rem1}),  specifically \( \mathcal{K}_3 \), that is a cycle of length 3. By definition of the $girth$ of a graph it follows that \( girth(\CC_f(\A)) = 3 \).
\end{proof}
\begin{rem}
	In Example \ref{ex1}, $F_5, F_6$ and $F_7$ are maximal filters of that residuated lattice $\A$. We can see that they form $\mathcal{K}_3$ as a subgraph of the graph $\CC_f(\A)$.
\end{rem}
\section{More on integer invariants of the graph $\CC_f(\A)$}
In this section, we aim to demonstrate that the graph \( \CC_f(\A) \) satisfies the condition where its clique number is equal to its chromatic number.

\begin{defn}\cite{B79}
A graph \( \Gamma \) is said to be weakly perfect if its chromatic number equals its clique number, that is, \( \chi(\Gamma) = \omega(\Gamma) \). 
\end{defn}
The next result parallels the classical Prime Avoidance Lemma for ideals in commutative rings \cite{K70}.
\begin{lem}\label{lem1}
Let $\A$ be a residuated lattice. Let \( F \subseteq A \) be a filter, and \( P_1, \dots, P_n \) be prime filters such that:
$ I \subseteq \bigcup_{i=1}^{n} P_i$. Then, there exists some \( j \in \{1, \dots, n\} \) such that: $F \subseteq P_j$.
\end{lem}
\begin{proof}
We proof the case of two prime filters, and it can be generalizes by induction.\\
Let \( F \) be a filter in \( \A \), and suppose $F \subseteq P_1 \cup P_2$ where \( P_1 \) and \( P_2 \) are prime filters.

We aim to show:
\[
F \subseteq P_1 \quad \text{or} \quad F \subseteq P_2
\]
By absurdum, suppose \( F \nsubseteq P_1 \) and \( F \nsubseteq P_2 \).  Then there exist elements \( x \in F \setminus P_1 \) and \( y \in F \setminus P_2 \).\\
Since \( F \) is a filter, it is upward closed and closed under finite meets, hence $x \wedge y \in F$. And we have \( x \notin P_1 \Rightarrow x \wedge y \notin P_1 \), and \( y \notin P_2 \Rightarrow x \wedge y \notin P_2 \). Therefore $
x \wedge y \notin P_1 \cup P_2 \Rightarrow a \wedge y \notin \bigcup_{i=1}^{2} P_i$ . 
Contradiction.\\
Hence, we conclude:
\[
F \subseteq P_1 \quad \text{or} \quad F \subseteq P_2
\]
The general case for \( n > 2 \) follows by induction.	
\end{proof}

\begin{thm} \label{thm1}
Let $\A$ be a finite residuated lattice. The comaximal filter graph \( \CC_f(\A) \) is a weakly perfect graph.
\end{thm}
\begin{proof}
To prove that for the graph $\CC_f(\A)$, the chromatic number is equal to the clique number, we are going to point out that they are both equal the number of maximal filters of $\A$.\\	
The definitions of  $\chi(\CC_f(\A))$ and  $\omega(\CC_f(\A))$ suggest trivialy that $\chi(\CC_f(\A))\geq \omega(\CC_f(\A))$. Assume that \( |Max(\A)| = n  \), and let \( Max(\A) = \{F_1, \ldots, F_n\} \). Since maximal filters are comaximal, consider the subgraph induced by \( \{F_1, \ldots, F_n\} \); it is the complete graph \( \mathcal{K}_n \), hence \(\omega(\CC_f(\A))\geq |Max(\A)| \), and we obtain $\chi(\CC_f(\A))\geq \omega(\CC_f(\A))\geq |Max(\A)|$.\\
What remains to be shown is that \( Max(\CC_f(\A))| \geq \chi(\CC_f(\A)) \).\\
Let us propose a coloring for graph $\CC_f(\A)$. By Lemma \ref{lem1}, the sets defined as follows:\\  
$
S_1 = \{ F \in V(\CC_f(\A)) \mid F \subseteq M_1 \}, 
S_2 = \{ F \in V(\CC_f(\A)) \mid F \subseteq M_2,\ F \notin S_1 \}, 
S_3 = \{ F \in V(\CC_f(\A)) \mid F \subseteq M_3,\ F \notin S_1 \cup S_2 \}, 
\ldots, 
S_n = \{ F \in V(\CC_f(\A)) \mid F \subseteq M_n,\ F \notin S_1 \cup \cdots \cup S_{n-1} \}
$\\
yields a proper \( n \)-coloring of the graph \( \CC_f(\A)  \), since every pair of distinct vertices in $S_i$ is non-adjacent. Therefore, \( |Max(\A)| = n \geq \chi(\CC_f(\A)) \), and it follows that $\chi(\CC_f(\A))$ and $\omega(\CC_f(\A))$ are equal to $|Max(\A)|$. This concludes the proof.
\end{proof}

\begin{pro}
Let $\A$ be a residuated lattice and $s$ be a positive integer greater than 1. Then the following assertions hold:
\begin{enumerate}
	\item Suppose \( |Max(\A)| = s \). Then the graph \( \CC_f(\A) \) admits a partition into \( s \) independent sets; that is, it is \( s \)-partite.
	\item If the graph \( \CC_f(\A) \) is \( s \)-partite, then \( |\text{Max}(R)| \leq s \).  
	Moreover, if \( \CC_f(\A) \) is not \( (s - 1) \)-partite, then \( |Max(\A)| = s \).
\end{enumerate}	
\end{pro}
\begin{proof}
(1) Assume that \( \text{Max}(\A) = \{ M_1, \ldots, M_s \} \), and define the sets $S_i$ as in the Proof of Theorem \ref{thm1}. According to Lemma \ref{lem1}, we have \( S_i \neq \emptyset \) for each \( i \). 	It is straightforward to verify that no two elements within any given \( S_i \) are adjacent in the graph. Which define a $s$-partite of the graph $\CC_f(\A)$.\\
(2) Let \( V_1, \ldots, V_s \) denote the \( s \) partitions of the vertex set of \(\CC_f(\A)\). Assume, for the sake of contradiction, that \( |Max(\A)| > n \).   Choose \( M_1, \ldots, M_{s+1} \in Max(\A) \). For each \( i =0...s+1 \), select an element $F_i \in S_i$ (where $S_i$ are defined as in the proof of Theorem \ref{thm1}). It is then easy to observe that the set \( \{F_1, \ldots, F_{s+1}\} \) forms a clique in \( \CC_f(\A) \). By the Pigeonhole Principle, at least two of the \( F_i \) must lie in the same partition \( V_j \), (where $1\leq j\leq s$) which contradicts the fact that no two vertices within a part $V_j$ are adjacent.  
Thus, we conclude that \( Max(\A)| \leq s \).
	
Now, suppose that \( \CC_f(\A) \) is not \( (s-1) \)-partite, and yet \( Max(\A)| = m < s \). From part (1), the graph must be \( m \)-partite, which again leads to a contradiction.
\end{proof}

\section{Some remarks on graph $\CC_f(\A)$, where $|\A|\leq 10$}
In \cite{BV10, K}, computer algorithms were employed to enumerate all non-isomorphic finite residuated lattices of size \( n \leq 12 \). We briefly describe all comaximal filter graphs constructed on specific residuated lattices (MV-algebra, G\"{o}del algebra, Heyting algebra, BL-algebra, and MTL-algebra) and provide some remarks for these graphs $\CC_f(\A)$ with the size of $\A$ less than or equal to 10. Since this investigation is carried out on residuated lattices generated by the algorithm in \cite{K}, an isomorphism between comaximal filter graphs does not necessarily imply an isomorphism between the underlying residuated lattices. Additionally, we give summary statistics of the obtained comaximal filters graphs. The results presented here are derived from constructions performed using algorithm \cite{A25}, as well as from some of the results in Section \ref{Sect2}, taking into account the constructions presented in \cite{BV10}.

\begin{pro}
	Let $\mathcal{A}$ be a residuated lattice of size 2, 3. Then $\CC_f(\A)$ is the null graph. 
\end{pro}
\begin{proof}
	By Proposition \ref{prop1} (1).
\end{proof}

\begin{pro}
Let $\mathcal{A}$ be a residuated lattice of size 4. There exists only one graph $\CC_f(\A)$ corresponding to an algebra $\A$ that is non null, and this graph is isomorphic to $\mathcal{K}_2$.
\end{pro}

\begin{pro}
Let $\A$ be an MV-algebra, a G\"{o}del algebra, a BL-algebra, or an MTL-algebra of size 5. Then $\CC_f(\A)$ is the null graph. 
\end{pro}
\begin{pro}
	Let $\A$ be a Heyting algebra or a non-prelinear residuated lattice of size 5. There exists only one graph $\CC_f(\A)$ for an algebra $\A$ that is non-null, and this graph is isomorphic to $\mathcal{K}_2$. 
\end{pro}

\begin{pro}
Let $\A$ be an MV-algebra, a Heyting algebra, or a non-prelinear residuated lattice of size 6. There exists only one graph $\CC_f(\A)$ for an algebra $\A$ that is non-null, and this graph is isomorphic to $\mathcal{K}_2$. 
\end{pro}
\begin{pro}\label{prop2}
	Let $\A$ be a Heyting algebra of size 6 not a G\"{o}del algebra. There exists only one graph $\CC_f(\A)$ for an algebra $\A$ that is  non-null, and this graph is isomorphic to $\mathcal{K}_2$.
	
\end{pro}
\begin{pro}\label{prop3}
Let $\A$ be a G\"{o}del-algebra of size 6. There exists only one graph $\CC_f(\A)$ for an algebra $\A$ that is non-null, and this graph is isomorphic to $\mathcal{P}_3$
	
\end{pro}
\begin{pro}\label{prop4}
Let $\A$ be a BL-algebra or an MTL-algebra of size 6. There exist two graphs $\CC_f(\A)$ corresponding to two algebras $\A$ that are non-null graphs. One is a Heyting algebra as stated in Proposition \ref{prop2} and the other is a G\"{o}del algebra as stated in Propositon \ref{prop3}.
\end{pro}
\begin{pro}
	Let $\A$ be a residuated lattice of size 6 that is not an MTL-algebra. There exist four graphs $\CC_f(\A)$ corresponding to algebras whose graphs are non-null. These graphs are all isomorphic to $\mathcal{K}_2$.
\end{pro}

\begin{pro}
	Let $\A$ be an MV-algebra or a G\"{o}del algebra of size 7. Then, $\CC_f(\A)$ is the null graph. 
\end{pro}
\begin{pro}
Let \(\mathcal{A}\) be a Heyting algebra of size 7. There exist three graphs \(\mathcal{C}_f(\mathcal{A})\) corresponding to algebras \(\mathcal{A}\) whose graphs are non-null. Two of them are isomorphic to \(\mathcal{K}_2\), and the third is isomorphic to \(\mathcal{P}_3\).
\end{pro}
\begin{pro}
	Let \(\mathcal{A}\) be a non-prelinear residuated lattice of size 7. There exist twenty one graphs \(\mathcal{C}_f(\mathcal{A})\) corresponding to algebras \(\mathcal{A}\) whose graphs are non-null. Each of these algebras has exactly two maximal filters. Twenty of them are isomorphic to \(\mathcal{K}_2\), and the remaining one is isomorphic to \(\mathcal{P}_3\).

\end{pro}

\begin{pro}
Let \(\mathcal{A}\) be an MV-algebra of size 8. There exist two graphs \(\mathcal{C}_f(\mathcal{A})\) corresponding to algebras \(\mathcal{A}\) whose graphs are non-null. One of them is isomorphic to \(\mathcal{K}_2\), and the second is given in Example \ref{ex1}.

\end{pro}
\begin{pro}
	Let \( A = \{0, a, b, c, d, 1\} \), and let \(\mathcal{A}\) be a Gödel algebra of size 8 defined on \( A \). There exist two graphs \(\mathcal{C}_f(\mathcal{A})\) corresponding to such algebras \(\mathcal{A}\) whose graphs are non-null.
	For one of them, the edge set is given by  
	$
	E(\mathcal{C}_f(\mathcal{A})) = \{ 
	\{f, 1\} - \{b, c, e, 1\},\ 
	\{b, c, e, 1\} - \{d, e, f, 1\},\ 
	\{b, c, e, 1\} - \{a, c, d, e, f, 1\} 
	\}.
	$	
	The second graph has three maximal filters, and its edge set is given by  
	$
	E(\mathcal{C}_f(\mathcal{A})) = \{
	\{c, 1\} - \{d, e, f, 1\},\ 
	\{e, 1\} - \{b, c, f, 1\},\ 
	\{f, 1\} - \{a, c, e, 1\},\ 
	\{f, 1\} - \{b, c, f, 1\},\ 
	\{a, c, e, 1\} - \{d, e, f, 1\},\ 
	\{b, c, f, 1\} - \{d, e, f, 1\} 
	\}.
	$

\end{pro}
\begin{pro}
Let \(\mathcal{A}\) be a Heyting algebra of size 8. There exist seven graphs \(\mathcal{C}_f(\mathcal{A})\) corresponding to algebras \(\mathcal{A}\) whose graphs are non-null. Three of them are isomorphic to \(\mathcal{K}_2\), two are isomorphic to \(\mathcal{P}_3\), and the last one is the graph given in Example \ref{ex1}.
	
\end{pro}
\begin{pro}
Let \(\mathcal{A}\) be a BL-algebra of size 8. There exist five graphs \(\mathcal{C}_f(\mathcal{A})\) corresponding to algebras \(\mathcal{A}\) whose graphs are non-null, but only two of them are non-idempotent and non-involutive. These last two are isomorphic to \(\mathcal{P}_3\).

\end{pro}
\begin{pro}
	Let \(\mathcal{A}\) be an MTL-algebra of size 8. There exist seven graphs \(\mathcal{C}_f(\mathcal{A})\) corresponding to algebras \(\mathcal{A}\) whose graphs are non-null, but only one of them is non-divisible, non-idempotent, and non-involutive. This last one is isomorphic to \(\mathcal{K}_2\).

\end{pro}
\begin{pro}

Let \(\mathcal{A}\) be a non-prelinear residuated lattice of size 8. There exist one hundred and fifteen graphs \(\mathcal{C}_f(\mathcal{A})\) corresponding to algebras \(\mathcal{A}\) whose graphs are non-null, but only one hundred and seven of them are isomorphic to \(\mathcal{K}_2\), and eight are isomorphic to \(\mathcal{P}_3\).

\end{pro}

\begin{pro}
Let $\A$ be a MV-algebra  of size 9. There exists one graph $\CC_f(\A)$ corresponding to  an algebra $\A$ isomorphic to $\mathcal{K}_2$.

\end{pro}
\begin{pro}
Let \(\mathcal{A}\) be a G\"{o}del algebra of size 9. There exists one graph \(\mathcal{C}_f(\mathcal{A})\) corresponding to an algebra \(\mathcal{A}\) whose graph is non-null. See Example \ref{ex3} for this graph.

\end{pro}
\begin{pro}
	Considering $A=\{0,a,b,c,d,e,f,g,1\}$ and $\A$ be a Heyting-algebra of size 9. There exist eleven graphs $\CC_f(\A)$ corresponding to algebras $\A$ whose graphs are non-null. One is given in Remark \ref{Rem2}, another one is given in Example \ref{ex3} and five are isomorphic to $\mathcal{K}_2$, three are isomorphic to $\mathcal{P}_3$ and the edges set of the last one is given by $E(\CC_f(\A))=\{\{f,g,1\}-\{b,c,e,g,1\}-\{b,e,f,g,1\}, \{b,c,e,g,1\}-\{a,c,d,e,f,g,1\}\}$.

\end{pro}
\begin{pro}
Let \(\mathcal{A}\) be a BL-algebra or MTL-algebra of size 9. There exist three graphs \(\mathcal{C}_f(\mathcal{A})\) corresponding to algebras \(\mathcal{A}\) whose graphs are non-null. Two are isomorphic to \(\mathcal{K}_2\), and the other one is given in Example \ref{ex3}.
\end{pro}

\begin{pro}
Let \(\mathcal{A}\) be an MV-algebra of size 10. There does not exist any \(\mathcal{C}_f(\mathcal{A})\) for any algebra \(\mathcal{A}\) that is a non-null graph.
\end{pro}
\begin{pro}
	Considering \( A = \{0, a, b, c, d, e, f, g, h, 1\} \) and let \(\mathcal{A}\) be a G\"{o}del algebra of size 10 on \(A\). There exist two graphs \(\mathcal{C}_f(\mathcal{A})\) corresponding to algebras \(\mathcal{A}\) whose graphs are non-null. The edge sets are given by:
	 $E(\CC_f(\A))=\{\{g,1\}-\{b,c,e,1\}-\{f,g,1\}, \{d,e,g,1\}-\{b,c,e,1\}-\{a,c,d,e,f,g,1\}\}$ and $E(\CC_f(\A))=\{\{f,g,1\}-\{b,c,e,g,1\}-\{d,e,f,g,1\}, \{b,c,e,g,1\}-\{a,c,d,e,f,g,1\}\}$.
\end{pro}
\begin{pro}
	Considering \( A = \{0, a, b, c, d, e, f, g, h, 1\} \) and let \(\mathcal{A}\) be a Heyting algebra of size 10 on \(A\). There are twenty-one graphs \(\mathcal{C}_f(\mathcal{A})\) corresponding to algebras \(\mathcal{A}\) whose graphs are non-null. Eight of them are isomorphic to \(\mathcal{K}_2\), seven are chains with three elements, and the edge sets for the six others are given by:\\
	 $E(\CC_f(\A))=\{\{d,e,g,1\}-\{b,c,e,1\}-\{f,g,1\}, \{g,1\}-\{b,c,e,1\}-\{a,c,d,e,f,g,1\}\}$,\\$E(\CC_f(\A))=\{\{f,g,1\}-\{b,c,e,g,1\}-\{d,e,f,g,1\}, \{b,c,e,g,1\}-\{a,c,d,e,f,g,1\}\}$,\\	 $E(\CC_f(\A))=\{\{g,1\}-\{b,c,e,1\}-\{a,c,d,e,f,g,1\}, \{d,e,g,1\}-\{b,c,e,1\}-\{f,g,1\}\}$,\\
	$E(\CC_f(\A))=\{\{f,g,1\}-\{b,c,e,g,1\}-\{d,e,f,g,1\}, \{b,c,e,g,1\}-\{a,c,d,e,f,g,1\}\}$,\\
	$E(\CC_f(\A))=\{\{d,e,f,g,1\}-\{c,g,1\}-\{e,g,1\}-\{b,c,f,g,1\}-\{d,e,f,g,1\}-\{a,c,e,g,1\}-\{c,g,1\}, \{a,c,e,g,1\}-\{b,c,f,g,1\}, \{a,c,e,g,1\}-\{e,g,1\}\}$,\\
	$E(\CC_f(\A))=\{\{f,1\}-\{b,c,d,g,1\}-\{a,c,d,f,1\}-\{g,1\}, \{b,c,d,g,1\}-\{e,f,g,1\}-\{a,c,d,f,1\}\}$, and\\ $E(\CC_f(\A))=\{\{d,e,1\}-\{f,g,1\}-\{a,c,d,e,g,1\}-\{b,c,e,f,g,1\}-\{d,e,1\}\}$.
\end{pro}
\begin{pro}
Considering \( A = \{0, a, b, c, d, e, f, g, h, 1\} \) and let \(\mathcal{A}\) be a BL-algebra of size 10 on \(A\). There are eight graphs \(\mathcal{C}_f(\mathcal{A})\) corresponding to algebras \(\mathcal{A}\) whose graphs are non-null. Two of them are isomorphic to \(\mathcal{K}_2\), three are isomorphic to $\mathcal{P}_3$ and the edge sets of the others are given by:\\
$E(\CC_f(\A))=\{\{g,1\}-\{b,c,e,1\}-\{a,c,d,e,f,g,1\}, \{d,e,g,1\}-\{b,c,e,1\}-\{f,g,1\}\}$,\\ 
	and\\
	 $E(\CC_f(\A))=\{\{f,g,1\}-\{b,c,e,g,1\}-\{d,e,f,g,1\}, \{b,c,e,g,1\}-\{a,c,d,e,f,g,1\}\}$. 
\end{pro}
\begin{pro}
	Considering \( A = \{0, a, b, c, d, e, f, g, h, 1\} \) and let \(\mathcal{A}\) be an MTL-algebra of size 10. There are seventeen graphs \(\mathcal{C}_f(\mathcal{A})\) corresponding to algebras \(\mathcal{A}\) whose graphs are non-null. Nine of them are isomorphic to \(\mathcal{K}_2\). The edge sets of the others are as follows:\\	
	- One is given by \( E(\mathcal{C}_f(\mathcal{A})) = \{\{f, g, 1\} - \{b, c, e, g, 1\} - \{d, e, f, g, 1\}, \{b, c, e, g, 1\} - \{a, c, d, e, f, g, 1\} \} \),\\
	- Six of them are isomorphic to $\mathcal{P}_3$.\\
	- The last one by \( E(\mathcal{C}_f(\mathcal{A})) = \{\{g, 1\} - \{b, c, e, 1\} - \{d, e, g, 1\}, \{f, g, 1\} - \{b, c, e, 1\} - \{a, c, d, e, f, g, 1\} \} \).
\end{pro}
\begin{thm}\cite{K30}
	A graph is planar if and only if it contains neither $\mathcal{K}_5$ or $\mathcal{K}_{3,3}$ as minors.
\end{thm}
\begin{rem}
	All the comaximal filter graphs constructed on residuated lattices of size less than or equal to 10 are planar.
\end{rem}

The first of the two following tables, Table \ref{Tab1}, taken from \cite{BV10}, which provides statistics on certain specific residuated lattices of size less than or equal to 10. Routine calculations allowed us to obtain Propositions 5.1–5.27, which are also confirmed by the algorithms developed by the authors \cite{A25}. All of this resulted to the creation of Table \ref{Tab2} which provides statistics on non-null comaximal filter graphs for some specific classes of residuated lattices of size less than or equal to 10. 

\begin{table}[H]
	\centering
\begin{tabular}{|c|c|c|c|c|c|c|c|c|c|c|}
	\hline 
	\rule[-1ex]{0pt}{2.5ex}  & 1 & 2 & 3 & 4 & 5 & 6 & 7 & 8 & 9 & 10 \\ 
	\hline 
	\rule[-1ex]{0pt}{2.5ex} All & 1 & 1 & 2 & 7 & 26 & 129 & 723 & 4712 & 34698 & 290565 \\ 
	\hline 
	\rule[-1ex]{0pt}{2.5ex} MTL & 1 & 1 & 2 & 7 & 23 & 99 & 464 & 2453 & 14087 & 88188 \\ 
	\hline 
	\rule[-1ex]{0pt}{2.5ex} BL & 1 & 1 & 2 & 5 & 9 & 20 & 38 & 81 & 160 & 326 \\ 
	\hline 
	\rule[-1ex]{0pt}{2.5ex} Heyting & 1 & 1 & 1 & 2 & 3 & 5 & 8 & 15 & 26 & 47 \\ 
	\hline 
	\rule[-1ex]{0pt}{2.5ex} G\"{o}del & 1 & 1 & 1 & 2 & 2 & 3 & 3 & 5 & 6 & 8 \\ 
	\hline 
	\rule[-1ex]{0pt}{2.5ex} MV & 1 & 1 & 1 & 2 & 1 & 2 & 1 & 3 & 2 & 2 \\ 
	\hline 
\end{tabular}
\caption{Number of residuated lattices of size $\leq$ 10}\label{Tab1} 
\end{table}

\begin{table}[H]
	\centering
	\begin{tabular}{|c|c|c|c|c|c|c|c|c|c|c|}
		\hline 
		\rule[-1ex]{0pt}{2.5ex}  & 1 & 2 & 3 & 4 & 5 & 6 & 7 & 8 & 9 & 10 \\ 
		\hline 
		\rule[-1ex]{0pt}{2.5ex} All & 0 & 0 & 0 & 1 & 1 & 6 & 21 & 122 & 693 & 4534 \\ 
		\hline 
		\rule[-1ex]{0pt}{2.5ex} MTL & 0 & 0 & 0 & 1 & 0 & 2 & 0 & 7 & 3 & 17 \\ 
		\hline 
		\rule[-1ex]{0pt}{2.5ex} BL & 0 & 0 & 0 & 1 & 0 & 2 & 0 & 5 & 3 & 8 \\ 
		\hline 
		\rule[-1ex]{0pt}{2.5ex} Heyting & 0 & 0 & 0 & 1 & 1 & 2 & 3 & 7 & 11 & 21 \\ 
		\hline 
		\rule[-1ex]{0pt}{2.5ex} G\"{o}del & 0 & 0 & 0 & 1 & 0 & 1 & 0 & 2 & 1 & 2 \\ 
		\hline 
		\rule[-1ex]{0pt}{2.5ex} MV & 0 & 0 & 0 & 1 & 0 & 1 & 0 & 2 & 1 & 0 \\ 
		\hline 
	\end{tabular}
	\caption{Number of non-null comaximal filter arising from residuated lattices of size $\leq$ 10}\label{Tab2} 
\end{table}


Graph theory in connection with algebraic structures has traditionally been developed around groups, semirings, and rings. More recently, researchers have extended this perspective to encompass broader algebraic structures, such as MV-algebras and BL-algebras. The present work is situated within this newer approach. Nevertheless, this does not preclude a possible return to specific classes of rings. Indeed, in papers such as \cite{BD09, BV10, BDM10, HLNN16, MASH21}, the authors characterize rings whose sets of ideals form specific types of residuated lattices. This deep connection between residuated lattices and ring theory is also central to the work presented in \cite{GY20} on zero-divisor graphs.

\end{document}